# The Maximal Integral Domain Generated By A Commutative Ring

Kerry M. Soileau

July 19, 2006

ksoileau@yahoo.com

ABSTRACT
In this paper, we exhibit the creation of the maximal integral domain $mid(R)$ generated by a nonzero commutative ring $R$.

Let $R$ be a nonzero commutative ring whose zero-divisors generate a <u>proper</u> ideal of $R$. $\mathbb{Z}_p$ is clearly such a one for prime $p$, whereas $\mathbb{Z}_6$ is not. Let $z(R)$ be the smallest prime ideal containing the zero-divisors of $R$. $z(R)$ exists because the maximal ideal containing the zero-divisors of $R$ exists and (like all maximal ideals) is prime. Thus the smallest prime ideal containing the zero-divisors of $R$ must exist.

Definition: $mid(R) \equiv R/z(R)$.

<u>Proposition</u>: $mid(R)$ is an integral domain.

Proof: Fix $x+z(R), y+z(R) \in R/z(R)$ with $(x+z(R))(y+z(R)) = 0+z(R)$. If $x+z(R) = 0+z(R)$ or $y+z(R) = 0+z(R)$ then we are done. Now assume $x+z(R) \neq 0+z(R)$ and $y+z(R) \neq 0+z(R)$. Then $x \notin z(R)$ and $y \notin z(R)$. Since $z(R)$ is prime, it follows that $xy \notin z(R)$, hence $0+z(R) = (x+z(R))(y+z(R)) = xy+z(R) \neq 0+z(R)$ which is a contradiction, hence we must have either $x+z(R) = 0+z(R)$ or $y+z(R) = 0+z(R)$. Therefore $mid(R)$ is an integral domain.